\documentclass[11pt]{article}

\usepackage[a4paper,top=1in,bottom=1in,left=1in,right=1in,footskip=0.25in,includeheadfoot]{geometry}
\usepackage{fancyhdr, titling, caption, amsmath, amsthm, amsfonts, mathtools, parskip, threeparttable, booktabs}
\usepackage[symbol]{footmisc}

\title{On Arithmetic Progressions and a Proof of the Nonexistence of Magic Squares of Squares}
\author{Oscar Hill}
\date{October 2025}

\newcommand{\overlap}{\,\,\|\,\,}
\newcommand{\loverlap}{<\!\!\|\,\,}
\newcommand{\mtnote}[1]{\textsuperscript{\TPTtagStyle{#1}}}

\begin{document}
\pagestyle{fancy}
\fancyhead{}
\renewcommand{\headrulewidth}{0pt}
\fancyhead[C]{\fontsize{9}{12} \selectfont On Arithmetic Progressions and a Proof of the Nonexistence of Magic Squares of Squares}

\maketitle

\begin{abstract}
\noindent
We explore some of the properties of consecutive, equally-summed arithmetic progressions of odd numbers, particularly their offsets and sums, 
before using them to prove that no \(3\times3\) magic squares of distinct square integers exist.
\bigskip\hrule\bigskip
\end{abstract}

\section{Introduction}
A \textit{magic square} is a square grid (of side length \(n\in\mathbb{N}\), called its \textit{order}) of distinct positive integers, whose rows, 
columns, major diagonal and minor diagonal all add to the same sum, called its \textit{magic constant}. The earliest concrete example is a 
\(3\times3\) magic square from the Chinese text \textit{Ta Tai Li-chi}, compiled in the first century CE [1] (Figure 1). Magic squares have since been 
studied by many cultures, initially often for mystical or divination purposes, but they are now researched for their mathematical properties.

\begin{center}
\begin{tabular}{ c|c|c }
4 & 9 & 2 \\
\hline
3 & 5 & 7 \\
\hline
8 & 1 & 6
\end{tabular}
\captionof{table}{\textsc{Figure 1}: The magic square described in the \textit{Ta Tai Li-chi} with a magic constant of 15}
\end{center}

A natural question to ask is, \textit{Does a magic square exist whose elements are all integers raised to the same (integral) power, \(d\)?} In 1770, 
Euler was the first to construct a \(4\times4\) magic square of square integers (Figure 2), and Rome and Yamagishi have since 
proven [2] that for every power \(d\geq2\), there exists an integer \(n_o(d)\) such that for all \(n\geq n_o(d)\) such a magic 
square of order \(n\) must exist. In particular, for \(d=2\), they showed that \(n\geq4\).

\begin{center}
\begin{tabular}{ c|c|c|c }
\(68^2\) & \(29^2\) & \(41^2\) & \(37^2\) \\
\hline
\(17^2\) & \(31^2\) & \(79^2\) & \(32^2\) \\
\hline
\(59^2\) & \(28^2\) & \(23^2\) & \(61^2\) \\
\hline
\(11^2\) & \(77^2\) & \(8^2\) & \(49^2\)
\end{tabular}
\captionof{table}{\textsc{Figure 2}: Euler's magic square of squares with a magic constant of 8515}
\end{center}

However, this result does not imply the nonexistence of a \(3\times 3\) magic square of squares. There have been efforts to find 
such a square, and many semi-magic squares of squares have been found (a square is considered \textit{semi-magic} if it contains repeated entries, 
has more than one magic constant, or both). One famous example is the Parker square [3] (Figure 3). Other \(3\times3\) semi-magic squares of 
squares have been discovered [4], but no true magic squares have ever been found.

\begin{center}
\begin{tabular}{ c|c|c }
\(29^2\) & \(1^2\) & \(47^2\) \\
\hline
\(41^2\) & \(37^2\) & \(1^2\) \\
\hline
\(23^2\) & \(41^2\) & \(29^2\)
\end{tabular}
\captionof{table}{\textsc{Figure 3}: The Parker square with magic constants 3051 and 4107 (along the minor diagonal)}
\end{center}

It is a well-known fact that the differences between consecutive square integers form an arithmetic progression of odd integers. Arithmetic 
progressions (which we shall henceforth call \textit{APs}) are particularly well-suited to this problem of magic squares of squares. They have been 
used before (for example, by Pierrat, Thiriet and Zimmermann [5]), but we approach this problem by considering pairs of consecutive APs of odd 
integers, where each AP has an equal sum.

\textbf{Overview.} In Section 2, we explore a few of the properties of such pairs of APs, namely the relationship between their offsets and the sum 
of each AP. We show how they relate to \(3\times3\) magic squares of squares in Section 3, before showing that they do not allow for the existence
of such a magic square.

\textbf{Notation.} We shall use the following notation in this paper (note that, unless stated otherwise, all APs have a difference of 2 between terms 
and contain only integers).

\begin{threeparttable}
\begin{tabular}{ c|l }
    \(\mathbb{N}\) & The set of natural numbers \(\{1,2,3,\ldots\}\) excluding 0 \\
    \(\mathbb{N}_0\) & The set of natural numbers \(\{0,1,2,3,\ldots\}\) including 0 \\
    \(\mathbb{N}_E\) & The set of even natural numbers \(\{0,2,4,6,\ldots\}\) \\
    \(\mathbb{Z}\) & The set of integers \(\{\ldots,-2,-1,0,1,2,\ldots\}\) \\
    \(\mathbb{Q}^+\) & The set of positive rationals \(a/b\), \(a,b\in\mathbb{N}\) \\
    \(\mathbb{R}^+\) & The set of positive real numbers \\
    \(\mathbf{p}(\mathbf{p}^o,\mathbf{p}^n)\) & An AP with starting offset\mtnote{*} \(\mathbf{p}^o\in\mathbb{N}_E\) and length \(\mathbf{p}^n\in\mathbb{N}\) \\
    \(\Sigma\mathbf{p}\) & The sum of the AP \(\mathbf{p}\), where \(\Sigma\mathbf{p}\in\mathbb{N}\) \\
    \(\mathcal{P}\) & A pair of consecutive APs \((\mathbf{p}_1,\mathbf{p}_2)\) with equal sums (an \textit{AP pair}) \\
    \(\Sigma\mathcal{P}\) & The sum of either of the APs in \(\mathcal{P}\), where \(\Sigma\mathcal{P}\in\mathbb{N}\) \\
    \(\mathcal{P}^{n_1}\) & Equivalent to \(\mathbf{p}_1^n\) in \(\mathcal{P}\) \\
    \(\mathcal{P}^{n_2}\) & Equivalent to \(\mathbf{p}_2^n\) in \(\mathcal{P}\) \\
    \(\mathcal{P}^1\) & The starting offset of \(\mathbf{p}_1\) in \(\mathcal{P}\) \\
    \(\mathcal{P}^2\) & The starting offset of \(\mathbf{p}_2\) in \(\mathcal{P}\), where 
        \(\mathcal{P}^2\vcentcolon=\mathcal{P}^1+2\mathcal{P}^{n_1}\) \\
    \(\mathcal{P}^3\) & The final offset contained by \(\mathbf{p}_2\) in \(\mathcal{P}\), where 
        \(\mathcal{P}^3\vcentcolon=\mathcal{P}^2+2(\mathcal{P}^{n_2}-1)\) \\
    \(\mathcal{P}_1<\mathcal{P}_2\) & Logically equivalent to \(\mathcal{P}_1^3<\mathcal{P}_2^1\) \\
    \(\mathcal{P}_1\overlap\mathcal{P}_2\) & (Read: \(\mathcal{P}_1\) \textit{overlaps} \(\mathcal{P}_2\)) Logically equivalent to 
        \(\mathcal{P}_1^1\leq\mathcal{P}_2^3\land\mathcal{P}_2^1\leq\mathcal{P}_1^3\) \\
    \(\mathcal{P}_1\loverlap\mathcal{P}_2\) & Logically equivalent to \(\mathcal{P}_1^1\leq\mathcal{P}_2^3\) \\
\bottomrule\addlinespace[1ex]
\end{tabular}
\begin{tablenotes}\footnotesize
    \item[*] See equation (2)
\end{tablenotes}
\end{threeparttable}

\section{Pairs of Consecutive APs}
\subsection{A Brief Overview of APs}
The sum \(S\in\mathbb{N}\) of any AP \(\mathbf{p}\) is given by
\begin{equation}
    S=\frac{n}{2}[2a+(n-1)d],
\end{equation}
where \(a\in\mathbb{N}\) is its first element, \(d\in\mathbb{N}\) is the difference between consecutive terms and \(n\in\mathbb{N}\) is the number of terms.
Since we are only interested in odd numbers, this equation can be simplified by adopting the above notation, setting \(d=2\) and \(a=1\), and adding an offset 
\(\mathbf{p}^o\in\mathbb{N}_E\) to each term in the progression, so that
\begin{equation}
    \Sigma\mathbf{p}=\mathbf{p}^n(\mathbf{p}^n+\mathbf{p}^o)
\end{equation}
and
\begin{equation}
    \mathbf{p}^o=\frac{\Sigma\mathbf{p}}{\mathbf{p}^n}-\mathbf{p}^n.
\end{equation}

\textbf{Lemma 2.1.} \textit{Given any \(x\in\mathbb{N}_0\), there is exactly one \(\mathbf{p}^o\) corresponding to \(x^2\) 
(where \(\mathbf{p}\) contains the differences between terms in a sequence of consecutive squares starting at \(x^2\)) with \(\mathbf{p}^o=2x\).}

\begin{proof}
The difference between \(x^2\) and \((x+1)^2\) is \(2x+1\), which can be contained within \(\mathbf{p}\). 
By definition, the initial term of \(\mathbf{p}\) is \(a\vcentcolon=2x+1\). By equating equations (1) and (2) (where \(d=2\)), \(\mathbf{p}^o=a-1\).
Therefore, \(\mathbf{p}^o=2x\).

\end{proof}

\textbf{Observation 2.2.} \textit{The sum of an AP \(\mathbf{p}(O,m)\), comprised of the consecutive APs \(\mathbf{p}_A(O,n)\) and
\(\mathbf{p}_B(O+2n,m-n),m>n\), is given by \(\Sigma\mathbf{p}=\Sigma\mathbf{p}_A+\Sigma\mathbf{p}_B\).}

\subsection{Properties of an AP Pair}
Let \(\mathcal{P}\) be an AP pair. Then, from equation (3),
\begin{equation}
    \mathcal{P}^1=\frac{\Sigma\mathcal{P}}{\mathcal{P}^{n_1}}-\mathcal{P}^{n_1}
\end{equation}
and
\begin{equation}
    \mathcal{P}^2=\frac{\Sigma\mathcal{P}}{\mathcal{P}^{n_2}}-\mathcal{P}^{n_2}=\mathcal{P}^1+2\mathcal{P}^{n_1}=
    \frac{\Sigma\mathcal{P}}{\mathcal{P}^{n_1}}+\mathcal{P}^{n_1}.
\end{equation}
However, \(\mathcal{P}^{n_1}\) can be written in terms of \(\mathcal{P}^{n_2}\) by introducing \(\kappa\in\mathbb{Q}^+\), defined as 
\(\kappa\vcentcolon=\mathcal{P}^{n_1}/\mathcal{P}^{n_2}\). Clearly, 
since \(\mathcal{P}^1<\mathcal{P}^2\), equation (2) shows that \(0<\mathcal{P}^{n_2}<\mathcal{P}^{n_1},\) so \(\kappa>1\). 
We now use equation (5) to show that
\begin{equation}
    \Sigma\mathcal{P}=(\mathcal{P}^{n_2})^2\kappa\frac{\kappa+1}{\kappa-1}.
\end{equation}
Therefore,
\begin{equation}
    \mathcal{P}^1=\mathcal{P}^{n_2}\left(\frac{\kappa+1}{\kappa-1}-\kappa\right)
\end{equation}
and
\begin{equation}
    \mathcal{P}^2=\mathcal{P}^{n_2}\left(\kappa\frac{\kappa+1}{\kappa-1}-1\right).
\end{equation}

\subsection[An Expression for k]{An Expression for \(\kappa\)}

Consider the expression for the sum of an AP, given by equation (6):
\begin{equation}\tag{6}
    \Sigma\mathcal{P}=(\mathcal{P}^{n_2})^2\kappa\frac{\kappa+1}{\kappa-1}.
\end{equation}
This can be simplified by defining a variable \(\alpha\in\mathbb{Q}^+,\)
\begin{equation}
    \alpha\vcentcolon=\kappa\frac{\kappa+1}{\kappa-1},
\end{equation}
where \(\alpha=\alpha_n/\alpha_d\), \(\alpha_n,\alpha_d\in\mathbb{N}\) 
(WLOG, assume \(\alpha_n\) and \(\alpha_d\) to be relatively prime). This can be used to show that
\begin{equation}
    \kappa=\frac{\alpha_n-\alpha_d\pm\sqrt{(\alpha_n-3\alpha_d)^2-8\alpha_d^2}}{2\alpha_d}.
\end{equation}
However, as \(\kappa\in\mathbb{Q}^+\), \((\alpha_n-3\alpha_d)^2-8\alpha_d^2\) must be square. 
Clearly, then, \(8\alpha_d^2\) is the difference between two squares, thus an AP \(\mathbf{q}(Q,N)\) can be used. From equation (2),
\begin{equation}
    N(N+Q)=8\alpha_d^2,
\end{equation}
where \(Q\) is the offset corresponding to \((\alpha_n-3\alpha_d)^2-8\alpha_d^2\). 
Let \(Q'\in\mathbb{N},Q'\vcentcolon=Q+2N\) 
be the offset corresponding to \((\alpha_n-3\alpha_d)^2\). Lemma 2.1 can then be used to calculate \(\alpha_n-3\alpha_d\) from 
\(Q'\), thus
\begin{equation}
    \alpha_n=\frac{4\alpha_d^2}{N}+3\alpha_d+\frac{N}{2}.
\end{equation}
Now, solving for \(\kappa\),
\begin{equation}
    \kappa_+=\frac{4\alpha_d+N_+}{N_+}
\end{equation}
and
\begin{equation}
    \kappa_-=\frac{2\alpha_d+N_-}{2\alpha_d}.
\end{equation}

\section{From AP Pairs to Magic Squares}
\subsection{The Link Between APs and Magic Squares of Squares}
Consider a \(3\times3\) magic square of \textbf{distinct} \(x_{11}\ldots x_{33}\in\mathbb{N}_0\):

\begin{center}
\begin{tabular}{ c|c|c }
\(x_{11}\) & \(x_{12}\) & \(x_{13}\) \\
\hline
\(x_{21}\) & \(x_{22}\) & \(x_{23}\) \\
\hline
\(x_{31}\) & \(x_{32}\) & \(x_{33}\)
\end{tabular}
\captionof{table}{\textsc{Figure 4}: An arbitrary \(3\times3\) magic square of distinct integers \(x_{11}\ldots x_{33}\)}
\end{center}
Let \(D_1,D_2\in\mathbb{Z}\). Using the fact that all rows, columns and both diagonals through \(x_{22}\) sum to the same magic constant \(K\in\mathbb{N}\),
\begin{align}
\begin{split}
    x_{31}&=x_{23}+D_1\\
    x_{12}&=x_{31}+D_1\\
    x_{11}&=x_{12}+D_2\\
    x_{22}&=x_{11}+D_1\\
    x_{33}&=x_{22}+D_1\\
    x_{32}&=x_{33}+D_2\\
    x_{13}&=x_{32}+D_1\\
    x_{21}&=x_{13}+D_1.\\
\end{split}
\end{align}
Equations (15) hold for all magic squares, but now suppose that all \(x_{ij}\) are square. Therefore, there must exist three
AP pairs \(\left(\mathcal{P}_1,\mathcal{P}_2,\mathcal{P}_3\right)\), each with the sum \(D_1\). There must also exist a pair of two APs (\textbf{not} an AP pair)
\(\left(\mathbf{p}_X(\mathcal{P}_1^3+2,(\mathcal{P}_2^1-\mathcal{P}_1^3)/2),\mathbf{p}_Y(\mathcal{P}_2^3+2,(\mathcal{P}_3^1-\mathcal{P}_2^3)/2)\right)\),
each with the sum \(D_2\) and placed between the AP pairs in the following order:
\begin{equation}
    \underbracket{D_1,D_1}_{\mathcal{P}_1},\underbracket{D_2}_{\mathbf{p}_X},\underbracket{D_1,D_1}_{\mathcal{P}_2},
    \underbracket{D_2}_{\mathbf{p}_Y},\underbracket{D_1,D_1}_{\mathcal{P}_3}.
\end{equation}
We assume that \(D_1,D_2>0\). If it happens that \(x_{31}<x_{23}\), the symmetry of (16) allows equations (15) to be rewritten so that 
\(D_1>0\), whilst maintaining (16). Likewise for \(D_2\) if \(x_{11}<x_{12}\), but it must be noted that this implies that 
either \(\mathcal{P}_1\overlap\mathcal{P}_2\overlap\mathcal{P}_3\) or \(\mathcal{P}_3<\mathcal{P}_2<\mathcal{P}_1\), 
and so \(\left(\mathbf{p}_X,\mathbf{p}_Y\right)\rightarrow\left(\mathbf{p}_X'(\mathcal{P}_2^1,(\mathcal{P}_1^3+2-\mathcal{P}_2^1)/2),
\mathbf{p}_Y'(\mathcal{P}_3^1,(\mathcal{P}_2^3+2-\mathcal{P}_3^1)/2)\right)\).

\subsection{The Existence of Magic Squares of Squares}
In Section 2, we explored some of the properties of AP pairs, including their
sums and the relationships between each offset. We shall use these properties to prove Theorem 3.1.

\textbf{Theorem 3.1.} \textit{It is impossible to construct a \(3\times3\) magic square consisting solely of square integers.}

\begin{proof}
Suppose that it is possible to construct a \(3\times3\) magic square of square integers. Since every element must be distinct, there must be a sequence 
of APs whose sums satisfy (16). There must therefore exist three \textbf{distinct} AP pairs, \(\mathcal{P}_1\), \(\mathcal{P}_2\) and \(\mathcal{P}_3\), 
whose respective sums are, from equations (6) and (9),
\begin{align}
    \Sigma\mathcal{P}_1&=\alpha_1(\mathcal{P}_1^{n_2})^2,\\
    \Sigma\mathcal{P}_2&=\alpha_2(\mathcal{P}_2^{n_2})^2
\end{align}
and
\begin{equation}
    \Sigma\mathcal{P}_3=\alpha_3(\mathcal{P}_3^{n_2})^2,\\
\end{equation}
where, from (16),
\begin{equation}
    \Sigma\mathcal{P}_1=\Sigma\mathcal{P}_2=\Sigma\mathcal{P}_3=D_1.
\end{equation}

Furthermore, define \(\beta_1,\beta_2\in\mathbb{Q}^+\), such that
\begin{equation}
    \beta_1\vcentcolon=\frac{\mathcal{P}_1^{n_2}}{\mathcal{P}_2^{n_2}}=\sqrt{\frac{\alpha_2}{\alpha_1}}
\end{equation}
and
\begin{equation}
    \beta_2\vcentcolon=\frac{\mathcal{P}_2^{n_2}}{\mathcal{P}_3^{n_2}}=\sqrt{\frac{\alpha_3}{\alpha_2}}.
\end{equation}

Also, let \(\beta_1=\beta_{1n}/\beta_{1d}\) and \(\beta_2=\beta_{2n}/\beta_{2d}\), where \(\beta_{1n},\beta_{1d},\beta_{2n},\beta_{2d}\in\mathbb{R}^+\). 
In terms of \(\alpha_1\), \(\alpha_2\) and \(\alpha_3\),
\begin{align}
    \beta_{1n}&=\sqrt{\frac{\alpha_{2n}}{\alpha_{1n}}},\\
    \beta_{1d}&=\sqrt{\frac{\alpha_{2d}}{\alpha_{1d}}},\\
    \beta_{2n}&=\sqrt{\frac{\alpha_{3n}}{\alpha_{2n}}}
\end{align}
and
\begin{equation}
    \beta_{2d}=\sqrt{\frac{\alpha_{3d}}{\alpha_{2d}}}.
\end{equation}
From (16), we also have \(\mathbf{p}_A(\mathcal{P}_1^2,(\mathcal{P}_2^1-\mathcal{P}_1^2)/2)\) and \(\mathbf{p}_B(\mathcal{P}_2^2,(\mathcal{P}_3^1-\mathcal{P}_2^2)/2)\) 
whose respective sums are, from equation (2),\footnote[1]{From (16), \(\mathcal{P}_2^1<\mathcal{P}_1^2\iff\mathcal{P}_3^1<\mathcal{P}_2^2\). Therefore, 
if \(\mathcal{P}_2^1<\mathcal{P}_1^2\), then we have \(\mathbf{p}_A'(\mathcal{P}_2^1,(\mathcal{P}_1^2-\mathcal{P}_2^1)/2)\) and\\
\(\mathbf{p}_B'(\mathcal{P}_3^1,(\mathcal{P}_2^2-\mathcal{P}_3^1)/2)\), and so \(\Sigma\mathbf{p}_A'=\left[(\mathcal{P}_1^2)^2-(\mathcal{P}_2^1)^2\right]/4\) and 
\(\Sigma\mathbf{p}_B'=\left[(\mathcal{P}_2^2)^2-(\mathcal{P}_3^1)^2\right]/4\), which leads to an equivalent proof.}
\begin{equation}
    \Sigma\mathbf{p}_A=\frac{1}{4}\left[(\mathcal{P}_2^1)^2-(\mathcal{P}_1^2)^2\right]
\end{equation}
and
\begin{equation}
    \Sigma\mathbf{p}_B=\frac{1}{4}\left[(\mathcal{P}_3^1)^2-(\mathcal{P}_2^2)^2\right].
\end{equation}

\textbf{Lemma 3.2.} \textit{Given the definitions for \(\mathbf{p}_A\) and \(\mathbf{p}_B\), \(\Sigma\mathbf{p}_A=\Sigma\mathbf{p}_B\).}

\begin{proof}
    There are three possible cases: (a) \(\mathcal{P}_1<\mathcal{P}_2<\mathcal{P}_3\), 
    (b) \(\mathcal{P}_1\overlap\mathcal{P}_2\overlap\mathcal{P}_3\) such that \(\mathcal{P}_2^1\geq\mathcal{P}_1^2\) and 
    \(\mathcal{P}_3^1\geq\mathcal{P}_2^2\), and (c) \(\mathcal{P}_3\loverlap\mathcal{P}_2\loverlap\mathcal{P}_1\) such that 
    \(\mathcal{P}_2^1<\mathcal{P}_1^2\) and \(\mathcal{P}_3^1<\mathcal{P}_2^2\). (It is clear, from (16) and equation (20), that one of
    these cases will hold for every \(\pm D_2\).) We prove that \(\Sigma\mathbf{p}_A=\Sigma\mathbf{p}_B\) for each case in turn.
    \begin{itemize}
        \item[(a)] Consider the APs \(\mathbf{p}_1(\mathcal{P}_1^2,\mathcal{P}_1^{n_2})\), \(\mathbf{p}_2(\mathcal{P}_2^2,\mathcal{P}_2^{n_2})\), 
        \(\mathbf{p}_1'(\mathcal{P}_1^3+2,(\mathcal{P}_2^1-\mathcal{P}_1^3)/2)\) and \\
        \(\mathbf{p}_2'(\mathcal{P}_2^3+2,(\mathcal{P}_3^1-\mathcal{P}_2^3)/2)\).
        From (16), \(\Sigma\mathbf{p}_1=\Sigma\mathbf{p}_2\) and \(\Sigma\mathbf{p}_1'=\Sigma\mathbf{p}_2'\). Since we have the APs 
        \(\mathbf{p}_A(\mathcal{P}_1^2,(\mathcal{P}_2^1-\mathcal{P}_1^2)/2)\) and \(\mathbf{p}_B(\mathcal{P}_2^2,(\mathcal{P}_3^1-\mathcal{P}_2^2)/2)\) then,
        from Observation 2.2, \\\(\Sigma\mathbf{p}_A=\Sigma\mathbf{p}_1+\Sigma\mathbf{p}_1'=\Sigma\mathbf{p}_2+\Sigma\mathbf{p}_2'=\Sigma\mathbf{p}_B\).
        \item[(b)] Consider the APs \(\mathbf{p}_1(\mathcal{P}_1^2,\mathcal{P}_1^{n_2})\), \(\mathbf{p}_2(\mathcal{P}_2^2,\mathcal{P}_2^{n_2})\), 
        \(\mathbf{p}_1'(\mathcal{P}_2^1,(\mathcal{P}_1^3+2-\mathcal{P}_2^1)/2)\) and \\
        \(\mathbf{p}_2'(\mathcal{P}_3^1,(\mathcal{P}_2^3+2-\mathcal{P}_3^1)/2)\).
        From (16), \(\Sigma\mathbf{p}_1=\Sigma\mathbf{p}_2\) and \(\Sigma\mathbf{p}_1'=\Sigma\mathbf{p}_2'\). Since we have the APs 
        \(\mathbf{p}_A(\mathcal{P}_1^2,(\mathcal{P}_2^1-\mathcal{P}_1^2)/2)\) and \(\mathbf{p}_B(\mathcal{P}_2^2,(\mathcal{P}_3^1-\mathcal{P}_2^2)/2)\) then,
        from Observation 2.2, \\
        \(\Sigma\mathbf{p}_1=\Sigma\mathbf{p}_A+\Sigma\mathbf{p}_1'=\Sigma\mathbf{p}_B+\Sigma\mathbf{p}_2'=\Sigma\mathbf{p}_2\). Hence,
        \(\Sigma\mathbf{p}_A=\Sigma\mathbf{p}_B\).
        \item[(c)] Consider the APs \(\mathbf{p}_1(\mathcal{P}_1^2,\mathcal{P}_1^{n_2})\), \(\mathbf{p}_2(\mathcal{P}_2^2,\mathcal{P}_2^{n_2})\), 
        \(\mathbf{p}_1'(\mathcal{P}_2^1,(\mathcal{P}_1^3+2-\mathcal{P}_2^1)/2)\) and \\
        \(\mathbf{p}_2'(\mathcal{P}_3^1,(\mathcal{P}_2^3+2-\mathcal{P}_3^1)/2)\).
        From (16), \(\Sigma\mathbf{p}_1=\Sigma\mathbf{p}_2\) and \(\Sigma\mathbf{p}_1'=\Sigma\mathbf{p}_2'\). Since we have the APs 
        \(\mathbf{p}_A(\mathcal{P}_2^1,(\mathcal{P}_1^2-\mathcal{P}_2^1)/2)\) and \(\mathbf{p}_B(\mathcal{P}_3^1,(\mathcal{P}_2^2-\mathcal{P}_3^1)/2)\) then,
        from Observation 2.2, \\
        \(\Sigma\mathbf{p}_1'=\Sigma\mathbf{p}_A+\Sigma\mathbf{p}_1=\Sigma\mathbf{p}_B+\Sigma\mathbf{p}_2=\Sigma\mathbf{p}_2'\). Hence,
        \(\Sigma\mathbf{p}_A=\Sigma\mathbf{p}_B\).
    \end{itemize}
\end{proof}

Equations (7), (8), (13) (or (14)) and (21-28), along with Lemma 3.2, can be used to show that
\begin{equation}
\begin{split}
    N_3^2\beta_{2n}^2\beta_{2d}^2\left(8\beta_{1d}^4\alpha_{1d}^2-N_2^2\right)^2
    &-N_2^2\left(8\beta_{2d}^4\beta_{1d}^4\alpha_{1d}^2-N_3^2\right)^2\\
    =4N_2^2N_3^2\beta_{1n}^2\beta_{1d}^2\beta_{2n}^2\beta_{2d}^2(\alpha_{1n}-\alpha_{1d})^2
    &-4N_2^2N_3^2\beta_{2n}^2\beta_{2d}^2\left(\beta_{1n}^2\alpha_{1n}-\beta_{1d}^2\alpha_{1d}\right)^2.
\end{split}
\end{equation}
Consider the RHS first. Using equation (12),
\begin{equation}
\begin{split}
    \mathrm{RHS}=&4N_2^2N_3^2\beta_{2n}^2\beta_{2d}^2\left(\beta_{1d}^2-\beta_{1n}^2\right)
    \left[\frac{16}{N_1^2}\beta_{1n}^2\alpha_{1d}^4
    +\frac{24}{N_1}\beta_{1n}^2\alpha_{1d}^3\right.\\
    &\left.+\left(13\beta_{1n}^2-\beta_{1d}^2\right)\alpha_{1d}^2
    +3N_1\beta_{1n}^2\alpha_{1d}
    +\frac{N_1^2}{4}\beta_{1n}^2\right].
\end{split}
\end{equation}
Compare this with the LHS of equation (29). Clearly, this is a quadratic in \(\alpha_{1d}^2\) and, as such, the coefficients of 
odd powers of \(\alpha_{1d}\) are \(0\), which must also hold in the RHS. Obviously, \(N_1,N_2,N_3,\beta_{1n},\beta_{1d},\beta_{2n},\beta_{2d}>0\), 
as per the definitions in equations (11), (21) and (22), so \(\beta_{1d}^2-\beta_{1n}^2=0\). However, if this is the case, then clearly \(\beta_1=1\). 
Therefore, from equation (21), \(\mathcal{P}_1^{n_2}=\mathcal{P}_2^{n_2}\). Appreciating that 
\(\Sigma\mathcal{P}_1=\Sigma\mathcal{P}_2=\Sigma\mathcal{P}_3\) and \(\Sigma\mathbf{p}_A=\Sigma\mathbf{p}_B\), equation (3) shows that 
\(\mathcal{P}_1^2=\mathcal{P}_2^2\) and hence that \(\mathcal{P}_1^1=\mathcal{P}_2^1=\mathcal{P}_3^1\). Therefore, 
\(\mathcal{P}_1=\mathcal{P}_2=\mathcal{P}_3\), contradicting our initial premise. Hence, a \(3\times3\) magic square consisting solely of square integers 
cannot be constructed.

\renewcommand\qedsymbol{QED}
\end{proof}

\section{References}
\begin{itemize}
    \item[] 
    {
        [1] Cammann S 1960 The Evolution of Magic Squares in China, \textit{Journal of the American Oriental Society} \textbf{80}, no. 2, pp. 116-124
    }
    \item[]
    {
        [2] Rome N and Yamagishi S 2024 On the Existence of Magic Squares of Powers \\
        arXiv:2406.09364v2 \textbf{[math.NT]}
    }
    \item[]
    {
        [3] Parker M 2016 The Parker Square - Numberphile \textit{Numberphile} [interview by Brady Haran] 
        https://www.youtube.com/watch?v=aOT\_bG-vWyg
    }
    \item[]
    {
        [4] Parker M 2025 A Magic Square Breakthrough - Numberphile \textit{Numberphile} [interview by Brady Haran] 
        https://www.youtube.com/watch?v=stpiBy6gWOA
    }
    \item[]
    {
        [5] Pierrat P, Thiriet F and Zimmermann P 2015 Magic Squares of Squares \textit{Loria, University of Lorraine}
    }
\end{itemize}
\end{document}